\documentclass[a4paper,12pt]{article}
\def\beq{\begin{equation}}
\def\eeq{\end{equation}}
\def\bea{\begin{eqnarray}}
\def\eea{\end{eqnarray}}
\def\nn{\nonumber}

\def\F{{\cal F}}
\def\g{{\bf g}}
\begin{document}
\thispagestyle{empty}
\vspace*{3cm}
\begin{center}
{\bf \Large Extensions of Peripheric Extended Twists 

and 

Inhomogeneous Lie Algebras}

\vspace{1cm}
N. Aizawa

\bigskip

Department of Applied Mathematics,

\medskip
Osaka Women's University, 

\medskip
Sakai, Osaka 590-0035, Japan

\medskip
(e-mail: aizawa@appmath.osaka-wu.ac.jp)
\end{center}

\vfill

\begin{abstract}
  Simple extensions of peripheric extended twists, introduced recently by 
Lyakhovsky and Del Olmo, are presented. Explicit form 
of twisting elements are given and it is shown that the new twists as well as 
peripheric extended twists are suitable to deform inhomogeneous Lie algebras such 
as $ isu(n),\ iso(n)$, $ (1+n) $ dimensional Schr\"odinger algebras and 
Poincar\'e algebra.

\end{abstract}

\newpage
%
%
%
\section{Introduction}

  It is no doubt that inhomogeneous Lie algebras play a crucial role in 
physics. One can readily imagine the Poincar\'e algebra, the Galilei 
algebra, the Schr\"odinger algebra and so on. 
If we expect that quantizations of these Lie algebras 
also play an important role in theoretical physics, various possibilities 
of quantization should be studied. In this article, we consider quantizations 
by means of Drinfel'd twist \cite{twist}, that preserves the triangularity of 
Lie algebras, along the line of recent developments [2-8]. 
Quantization of a Lie algebra $ \g$ by Drinfel'd twist is defined by an 
invertible element $ \F \in U(\g) \otimes U(\g) $ subject to the relations
\bea
 & & \F_{12} (\Delta \otimes id)(\F) = \F_{23} (id \otimes \Delta)(\F), \label{defF} \\
 & & (\epsilon \otimes id)(\F) = (id \otimes \epsilon)(\F) = 1, \nn
\eea
where $ \Delta $ and $ \epsilon $ denote the coproduct and the counit of 
$ U(\g) $, respectively. The quantized Lie algebra $ U_{\F}(\g) $ has 
the same commutation relations as $\g$ and deformed coproducts 
$ \Delta_{\F} = \F \Delta \F^{-1}. $ 
Explicit forms of the twist elements $ \F $ are important because not only the 
coproducts but also deformation of all other quantities are caused by the 
twisting element. Especially, the universal $R$-matrix for $ U_{\F}(\g) $ 
is given by $ R_{\F} = \F_{21} \F^{-1}_. $ Note that irreducible representations 
for $ U_{\F}(\g) $ are exactly the same as $ \g $ because of undeformed 
commutation relations. Therefore once an explicit form of twist element is 
obtained, matrix solutions of quantum Yang-Baxter equation are immediately 
calculated, so that it is an easy exercise to construct a quantum group dual 
to $ U_{\F}(\g) $, differential calculi covariant under the quantum group 
and so on. Some physical application of twisting are discussed in \cite{KS,MO}.

  There are not so many literatures (to the author's knowledge) discussing twist 
deformation of inhomogeneous Lie algebras. Twist deformation of the 
Poincar\'e algebra based on an abelian subalgebra is discussed 
in \cite{twistedPoin}. In Ref.\cite{Mur}, the Jordanian twists \cite{jor1,jor2} 
are generalized to multidimension and application to the Poincar\'e algebra is 
considered. 
Other approaches to obtain triangular deformation of 
inhomogeneous Lie algebras are found in [15-19] (see also the
references therein). We construct, in this article, explicit forms of 
twisting elements that are applicable to various inhomogeneous Lie algebras 
including the Poincar\'e algebra. Our twisting elements are extensions of 
the peripheric extended twists (PET). The PET was introduced in \cite{pet} 
as nontrivial limits of extended Jordanian twists \cite{et}. 
Four dimensional subalgebras $ {\bf h} \subset \g $ are considered in both 
the PET and the extended Jordanian twists. Then, regarding $ \F \in U({\bf h}) \otimes U({\bf h})$ 
as a twisting element for $ U(\g) $, a twisted algebra $ U_{\F}(\g) $ 
is constructed by the $ \F.$ 

  The plan of this article is as follows: In the next section, 
after a brief review of PET, we extend the PET 
by adding one more generator to the four dimensional subalgebra of PET. 
As a result, we obtain four different twisting elements. In \S 3, the 
twisting elements are applied to $ isu(n),\ iso(n), $ Schr\"odinger algebras and 
Poincar\'e algebra. \S 4 is a conclusion. 

%
%
%
\section{Peripheric extended twists and their extensions}
\setcounter{equation}{0}
\subsection{Peripheric extended twists}

  The Jordanian twists for a Lie algebra ${\bf g}$ is a quantization 
using a Borel subalgebra $ \{H, E \} \in {\bf g} \ $\cite{jor1,jor2}. 
By adding two more elements 
$ A, B $ to the Borel subalgebra, the subalgebra used in extended Jordanian twists 
is obtained \cite{et}.  The extended Jordanian twists have two nontrivial limits specified by subalgebras 
$ {\bf L}^c $ and ${\bf L}'^c $ \cite{pet}. The elements of subalgebras ${\bf L}^c $ 
and $ {\bf L}'^c $ are denoted again by $ \{H, E, A, B \}. $ 
The algebra $ {\bf L}^c $ is defined by the commutation relations
\bea
 & &  [H,\; E] = \delta E, \qquad [H,\; A] = 0, \qquad [H,\; B] = \delta B, \nn \\
 & &  [A,\; B] = \gamma E, \qquad [E,\; A] = [E,\; B] = 0,    \label{Lc}
\eea
and $ {\bf L}'^c $ by 
\bea
 & & [H,\; E] = \delta E, \qquad [H,\; A] = \delta A, \qquad 
     [H,\; B] = 0, \nn \\
 & & [A,\; B] = \gamma E,\qquad [E,\; A] = [E,\; B] = 0, \label{Lcd}
\eea
where $ \gamma $ and $ \delta $ are arbitrary nonzero complex numbers. 
One can say that $ {\bf L}^c $ and $ {\bf L}'^c $ are isomorphic, since 
\beq
 H \leftrightarrow H, \qquad E \leftrightarrow E, \qquad 
 A \leftrightarrow -B, \qquad B \leftrightarrow A, \label{LcandLcd}
\eeq
give the isomorphism. However, it turns out that ${\bf L}^c$ and 
$ {\bf L}'^c $ give different twists as limits of extended Jordanian twists \cite{pet}. 
We thus treat these two cases separately. 

  In the case of $ {\bf L}^c, $ the twisting element of PET is given by
\beq
  \F_P = \Phi_P \Phi_j, \qquad \Phi_P = \exp(A\otimes B e^{-\delta \sigma}), 
  \qquad \Phi_j = e^{H \otimes \sigma},           \label{twistP}
\eeq
where 
\beq
  \sigma = \frac{1}{\delta} \ln(1 + \gamma E),               \label{sigma}
\eeq
and $ \Phi_j $ is the twisting element of the Jordanian twist. The twisting 
element $ \F_P $ leads to the algebra $U_P({\bf L}^c) $ with the twisted 
coproduct $ \Delta_P \equiv \F_P \Delta \F_P^{-1} $:
\bea
 & & \Delta_P(H) = H \otimes e^{-\delta \sigma} + 1 \otimes H 
     - \delta A \otimes Be^{-2\delta\sigma}, \nn \\
 & & \Delta_P(A) = A \otimes e^{-\delta\sigma} + 1 \otimes A, \label{coproLc} \\
 & & \Delta_P(B) = B \otimes e^{\delta\sigma} + e^{\delta\sigma} \otimes B, \nn \\
 & & \Delta_P(E) = E \otimes e^{\delta\sigma} + 1 \otimes E. \nn 
\eea
It is important to note that there exist two primitive elements in $ U_P({\bf L}^c) $, 
namely, $ \sigma $ and $ Be^{-\delta\sigma}_. $ 

  In the case of $ {\bf L}'^c $, the twisting element is given by
\beq
  \F_{P'} = \Phi_{P'} \Phi_j, \qquad \Phi_{P'} = e^{A \otimes B}, \label{twistPd}
\eeq
where $ \Phi_j $ is same as (\ref{twistP}). We denote the twisted algebra 
by $ U_{P'}({\bf L}'^c) $ and the twisted coproducts by 
$ \Delta_{P'} \equiv \F_{P'} \Delta \F_{P'}^{-1} $:
\bea
 & & \Delta_{P'}(H) = H \otimes e^{-\delta\sigma} + 1 \otimes H 
     - \delta A \otimes Be^{-\delta\sigma}, \nn \\
 & & \Delta_{P'}(A) = A \otimes 1 + 1 \otimes A, \label{coproLcd} \\
 & & \Delta_{P'}(B) = B \otimes 1 + e^{\delta\sigma} \otimes B, \nn \\
 & & \Delta_{P'}(E) = E \otimes e^{\delta\sigma} + 1 \otimes E. \nn
\eea
Again there exist two primitive elements in $ U_{P'}({\bf L}'^c) $, 
namely, $ \sigma $ and $ A. $ 

  Both twisting elements $ \F_{P} $ and $ \F_{P'} $ satisfy the relations
\beq
 (\Delta \otimes id)(\F) = \F_{13} \F_{23}, \qquad 
 (id \otimes \Delta_{\F})(\F) = \F_{12} \F_{13}.  \label{fac}
\eeq
These relations guarantee that the twisting elements satisfy (\ref{defF}). 
It is also possible to regard $ \Phi_{P} $ and $ \Phi_{P'} $ as twisting elements 
for the Jordanian deformed algebras, that is, we start with algebras with 
twisted coproduct $ \Delta_j = \Phi_j \Delta \Phi_j^{-1} $, then consider the 
twisting by $ \Phi_{P} $ or $ \Phi_{P'}. $ In such situation, it turns out that 
$ \Phi_{P} $ and $ \Phi_{P'} $ have different factorizable properties stems from 
the fact that $ U_{P}({\bf L}^c) $ and $ U_{P'}({\bf L}'^c) $ have different 
primitive elements. The twisting element $ \Phi_{P} $ satisfies the same relation 
as (\ref{fac}), while $ \Phi_{P'} $ satisfies
\beq
  (\Delta_{P'} \otimes id)(\Phi_{P'}) = (\Phi_{P'})_{13} (\Phi_{P'})_{23}, 
  \qquad
  (id \otimes \Delta_j)(\Phi_{P'}) = (\Phi_{P'})_{12} (\Phi_{P'})_{13}. \label{fac2}
\eeq
Although this is not mentioned in \cite{pet}, these PET are appropriate for 
twisting inhomogeneous Lie algebras as we shall see later.

\subsection{Extensions of PET}

  The existence of primitive elements except $ \sigma $ allows us to 
extend the PET with an extra Jordanian like factor. 
To this end, we need one more generator $J$ and 
consider five dimensional subalgebras: $ {\bf L} = {\bf L}^c \cup \{J\} $ and 
$ {\bf L}' = {\bf L}'^c \cup \{J\}. $ Additional commutation relations to 
define $ {\bf L} $ are given by
\beq
 [J,\; H] = [J,\; E] = 0, \qquad [J,\; A] = -\mu A, \qquad [J,\; B] = \mu B,
 \label{Ladd}
\eeq
and for $ {\bf L}' $
\beq
 [J,\; H] = [J,\; E] = 0, \qquad [J,\; A] = \mu A, \qquad [J,\; B] = -\mu B,
 \label{Ldadd}
\eeq
where $ \mu $ is an arbitrary complex constant. We define for $ {\bf L} $
\beq 
 \Phi = e^{J \otimes \rho}, \qquad \rho = \frac{1}{\mu}\ln(1+\mu Be^{-\delta\sigma}),
 \label{phi}
\eeq
and for $ {\bf L}' $
\beq
  \Phi' = e^{J \otimes \rho'}, \qquad
  \rho' = \frac{1}{\mu} \ln(1 + \mu A). \label{phid}
\eeq
Note that $\{J, Be^{-\delta\sigma} \} $ in (\ref{phi}) and $ \{J, A\} $ in 
(\ref{phid}) satisfy the commutation relation of Borel subalgebra, respectively. 
Then it can be verified that followings are twisting elements for $ {\bf L} $
\beq
 \F = \Phi \F_{P}, \qquad \tilde \F = \Phi_{21} \F_{P},   \label{twistorL}
\eeq
and for $ {\bf L}' $
\beq
 \F' = \Phi'\F_{P'}, \qquad \tilde \F' = \Phi'_{21} \F_{P'}. \label{twistorLd}
\eeq
To show this, we first calculate PET for $J$. It is easy to see that $J$ is primitive 
with respect to both twisting by $ {\bf L}^c$ and $ {\bf L}'^c$:
\beq
 \Delta_{P}(J) = \Delta_{P'}(J) = J \otimes 1 + 1 \otimes J. \label{coproJ}
\eeq
We next calculate twisted coproduct for the elements of ${\bf L}$ and ${\bf L}'$ 
by $ \F (\tilde \F) $ and $ \F' (\tilde \F') $, respectively. 
Twisted coproducts of $ {\bf L} $ by $ \F $ are given by
\bea
 & & \Delta_{\F}(J) = J \otimes e^{-\mu\rho} + 1 \otimes J, \nn \\
 & & \Delta_{\F}(H) = H \otimes e^{-\delta\sigma} + 1 \otimes H - J \otimes e^{-\mu\rho} 
     \{ (\delta-1) Be^{-\delta\sigma} + Be^{-2\delta\sigma} \} \nn \\
 & & \qquad \quad \; \, - \delta A \otimes Be^{-2\delta\sigma-\mu\rho}, \nn \\
 & & \Delta_{\F}(A) = A \otimes e^{-\delta\sigma-\mu\rho} + 1 \otimes A
     - \gamma J \otimes Ee^{-\delta\sigma - \mu\rho}, \label{copro1} \\
 & & \Delta_{\F}(B) = B \otimes e^{\delta\sigma+\mu\rho} + e^{\delta\sigma} \otimes B,
     \nn \\
 & & \Delta_{\F}(E) = E \otimes e^{\delta\sigma} + 1 \otimes E. \nn
\eea
Twisted coproducts of $ {\bf L} $ by $ \tilde\F$ are given by
\bea
 & & \Delta_{\tilde\F}(J) = J \otimes 1 + e^{-\mu\rho} \otimes J, \nn \\
 & & \Delta_{\tilde\F}(A) = A \otimes e^{-\delta\sigma} + e^{-\mu\rho} \otimes A 
     - \gamma Ee^{-\delta\sigma-\mu\rho} \otimes Je^{-\sigma\mu}, \label{copro2} \\
 & & \Delta_{\tilde\F}(B) = B \otimes e^{\delta\sigma} + e^{\delta\sigma+\mu\rho} 
     \otimes B, \nn \\
 & & \Delta_{\tilde\F}(E) = E \otimes e^{\delta\sigma} + 1 \otimes E. \nn
\eea
It seems to be difficult to have a closed form of $ \Delta_{\tilde\F}(H) $ 
because of the last term of $ \Delta_{P}(H). $ We therefore do not give it here. 
Twisted coproducts of $ {\bf L}' $ by $ \F' $ are given by
\bea
 & & \Delta_{\F'}(J) =  J \otimes e^{-\mu\rho'} + 1 \otimes J, \nn \\
 & & \Delta_{\F'}(A) = A \otimes e^{\mu\rho'} + 1 \otimes A, \label{copro3} \\
 & & \Delta_{\F'}(B) = B \otimes e^{-\mu\rho'} + e^{\delta\sigma} \otimes B 
     + \gamma Je^{\delta\sigma} \otimes Ee^{-\mu\rho'}, \nn \\
 & & \Delta_{\F'}(E) = E \otimes e^{\delta\sigma} + 1 \otimes E. \nn
\eea
It also seems to be difficult to have a closed form of $ \Delta_{\F'}(H) $ 
because of the last term of $ \Delta_{P'}(H). $ We do not give it here. 
Twisted coproducts of $ {\bf L}' $ by $ \tilde \F' $ are given by
\bea
 & & \Delta_{\tilde\F'}(J) = J \otimes 1 + e^{-\mu\rho'} \otimes J, \nn \\
 & & \Delta_{\tilde\F'}(H) = H \otimes e^{-\delta\sigma} + 1 \otimes H 
     + \frac{\delta}{\mu} (e^{-\mu\rho'}-1) \otimes Je^{-\delta\sigma} 
     - \delta A e^{-\mu\rho'} \otimes Be^{-\delta\sigma}, \nn \\
 & & \Delta_{\tilde\F'}(A) = A \otimes 1 + e^{\mu\rho'} \otimes A, \label{copro4} \\
 & & \Delta_{\tilde\F'}(B) = B \otimes 1 + e^{\delta\sigma -\mu\rho'}\otimes B 
     + \gamma Ee^{-\mu\rho'} \otimes J, \nn \\
 & & \Delta_{\tilde\F'}(E) = E \otimes e^{\delta\sigma} + 1 \otimes E. \nn
\eea
From these relations, it is seen that $ \sigma, \rho $ and $ \rho' $ are primitive 
in all cases. Together with (\ref{coproJ}), it follows that
\beq
  (\Delta_{\alpha} \otimes id)(\Psi) = \Psi_{13} \Psi_{23}, \qquad
  (id \otimes \Delta_f)(\Psi) = \Psi_{12} \Psi_{13}, \label{fact1}
\eeq
for $ (\alpha, f, \Psi) = (P, \F, \Phi), (P', \F', \Phi') $ and 
\beq
  (\Delta_f \otimes id)(\Psi) = \Psi_{13} \Psi_{23}, \qquad
  (id \otimes \Delta_{\alpha})(\Psi) = \Psi_{12} \Psi_{13}, \label{fact2}
\eeq
for $ (\alpha, f, \Psi) = (P, \tilde\F, \Phi_{21}), (P', \tilde\F', \Phi'_{21}). $ 
The relations (\ref{fact1}) and (\ref{fact2}) imply that $ \Phi $ and $ \Phi_{21} $ 
(resp. $ \Phi' $ and $ \Phi'_{21} $) are twisting elements for $ U_{P}({\bf L}) $ 
(resp. $ U_{P'}({\bf L}')). $ It is, indeed, easy to verify the second relation of 
(\ref{defF}). Thus the statement has been proved. 

  We can introduce a deformation parameter $z$ in the above expressions by 
the following replacement
\beq
   \left\{
   \begin{array}{lcl}
   E \rightarrow zE,\ B \rightarrow zB, & \hspace{5mm} & {\rm for} \ \ {\bf L}, \\
                    & & \\
   E \rightarrow zE,\ A \rightarrow zA, & & {\rm for} \ \ {\bf L}'.
   \end{array}
   \right.            \label{defpara}
\eeq
Undeformed algebras are recovered in the limit of $ z = 0. $ 
Let us consider the classical $r$-matrices obtained from the universal $R$-matrices 
given by twisting elements (\ref{twistorL}) and (\ref{twistorLd}). Classical 
$r$-matrices are obtained by 
keeping up to the first order in $z$ in the expansion 
of the universal $R$-matrices. The classical $r$-matrices obtained from 
$ {\cal R} = \F_{21} \F^{-1} $ and 
$ \tilde {\cal R} = \tilde \F_{21} \tilde \F^{-1} $ are 
\beq
  \pm J \wedge B + A \wedge B + \frac{\gamma}{\delta} H \wedge E, 
  \label{cr1}
\eeq
where $+ \;(-)$ corresponds to $ {\cal R} \; (\tilde {\cal R}). $ 
The classical $r$-matrices obtained from $ {\cal R}'= \F'_{21} \F'^{-1} $ and 
$ \tilde {\cal R}' = \tilde \F'_{21} \tilde \F'^{-1} $ are
\beq
  \pm J \wedge A + A \wedge B + \frac{\gamma}{\delta} H \wedge E, 
  \label{cr2}
\eeq
where $+\; (-)$ corresponds to $ {\cal R}'\; (\tilde {\cal R}'). $ 
These classical $r$-matrices (\ref{cr1}) and (\ref{cr2}) solve the classical 
Yang-Baxter equation.

%
%
%
\section{Twist deformation of inhomogeneous Lie algebras}
\setcounter{equation}{0}

  In this section, we show that many inhomogeneous Lie algebras can be 
twisted by the twisting elements (\ref{twistorL}) and (\ref{twistorLd}). 
The algebras considered are $ isu(n),\ iso(n), $ the Shor\"odinger algebra 
in $(1+n)$ spacetime. We also treat the Poincar\'e algebra in $(1+3)$ spacetime 
separately because of its physical importance. What is necessary is to show that 
these algebras have the subalgebras $ {\bf L} $ and $ {\bf L}'. $ It is, however, 
enough to show the existence of $ {\bf L} $, because the subalgebra 
$ {\bf L}' $ is obtained by making use of the replacement (\ref{LcandLcd}). 
Note that $ {\bf L} $ and ${\bf L}'$ contain the subalgebra $ {\bf L}^c $ and 
$ {\bf L}'^c $ of PET so that inhomogeneous Lie algebras considered in this section 
admit PET, too. 

  The Lie algebra $ isu(n) $ is defined by
\bea
 & & [U^a{}_b,\; U^c{}_d] = U^a{}_d \delta_b{}^c - U^c{}_b \delta_d{}^a, \nn \\
 & & [U^a{}_b,\; P^c] = P^a \delta_b{}^c, \qquad
     [U^a{}_b,\; P_c] = -P_b \delta^a{}_c, \label{iun} \\
 & & [P^a,\; P^b] = [P^a,\; P_b] = [P_a,\; P_b] = 0, \nn
\eea
where $ a, b, c, d \in \{1, 2, \cdots, n\}. $
The subalgebra $ {\bf L} $ exists for $ n \geq 4 $ and given by
\bea
 & & J = \sum_{2 \leq k \leq n/2}(U^k{}_{n-k} - U^{n-k}{}_k), \qquad 
     H = U^1{}_1-U^n{}_n, \qquad 
     E = P^1 + P_n, \nn \\
 & & A = \sum_{2 \leq k \leq n/2} \alpha_k \{P^k-P_k-i(P^{n-k}-P_{n-k})\}, \label{subisun} \\
 & & B = \sum_{2 \leq k \leq n/2} \beta_k \{U^1{}_k + U^k{}_n + i(U^1{}_{n-k}
       +U^{n-k}{}_n)\}, \nn
\eea
where $ \alpha_k $ and $ \beta_k $ are arbitrary constant satisfying the relation
\beq
  \gamma = -2 \sum_{2 \leq k \leq n/2} \alpha_k \beta_k.      \label{const1}
\eeq
Other parameters in $ {\bf L} $ are $ \delta = 1, \mu = i.$ 

  The Lie algebra $ iso(n) $ is generated by $ Y_{ab} = -Y_{ba} $ and $ P_a $ 
satisfying the relations
\bea
 & & [Y_{ab},\; Y_{cd}] = Y_{ad} \delta_{bc} + Y_{bc} \delta_{ad} 
     - Y_{ac} \delta_{bd} - Y_{bd} \delta_{ac}, \label{ison} \\
 & & [Y_{ab},\; P_c] = P_a \delta_{bc} - P_b \delta_{ac}, 
     \qquad
     [P_a,\; P_b] = 0, \nn
\eea
where $ a, b, c, d \in \{1, 2, \cdots, n\}. $ The subalgebra $ {\bf L} $ is found 
for $ n \geq 4 $ and given by
\bea
 & & J = \sum_{2 \leq k \leq n/2} Y_{k\; n-k+1}, \qquad 
     H = Y_{1n}, \qquad E = P_1 + iP_n, \nn \\
 & & A = \sum_{2 \leq k \leq n/2} \alpha_k (P_k -iP_{n-k+1}), \label{subison} \\
 & & B \sum_{2 \leq k \leq n/2} \beta_k (Y_{1k}-iY_{kn}+i Y_{1\; n-k+1} + Y_{n-k+1\; n}), \nn
\eea
where $ \alpha_k $ and $ \beta_k $ are arbitrary constant satisfying the same relation 
as (\ref{const1}). 
Other parameters in $ {\bf L} $ are $ \delta = \mu = i.$ 

  We next consider the centrally extended Schr\"odinger algebra in $(1+n)$ 
spacetime \cite{DDM}. It is generated by $ \frac{1}{2}n(n+3)+4 $ elements, namely, 
$P_t$: time translation, $P_a$: space translations, $G_a$: Galilei transformations, 
$ J_{ab}$: rotations, $K$: conformal transformation, $D$: dilatation and 
$M$: center (mass). All suffixes range from $1$ to $n$. The nonvanishing 
commutation relations are given by
\beq
\begin{tabular}{lll}
 $[P_t, D] = 2P_t, \quad $& $ [P_t, G_a] = P_a, \quad $&$ [P_t, K] = D,$ \\
 $[P_a, D] = P_a, \quad $&$ [P_a, K] = G_a, \quad $&$ [D, G_a] = G_a,$  \\
 $[D, K] = 2K, \quad $ &$ [P_a, G_b] = \delta_{ab} M,$ & \\
 $[P_a, J_{bc}] = \delta_{ac} P_b - \delta_{ab} P_c,$ & 
 $[G_a, J_{bc}] = \delta_{ac} G_b - \delta_{ab} G_c,$ & \\
 \multicolumn{2}{l}{$[J_{ab}, J_{cd}] = \delta_{ac}J_{bd} + \delta_{bd}J_{ac}
 -\delta_{ad}J_{bc} - \delta_{bc}J_{ad}.$} & 
\end{tabular}
\label{schdef}
\eeq
The subalgebra $ {\bf L} $ is found for $ n \geq 2 $:
\bea
 & & J = i \sum_{<k>} J_{k\; k+1} + D, \qquad H = i \sum_{<k>} J_{k\; k+1} + D, 
     \qquad B = P_t, \nn \\
 & & A = \sum_{<k>} \alpha_k (G_k + iG_{k+1}), \qquad 
     E = \sum_{<k>} \alpha_k (P_k + iP_{k+1}), \label{subsch}
\eea
where $\displaystyle{\sum_{<k>}} $ means that $k$ takes odd integers ranging 
from $1$ to $ n-1 $ (or $ n-2 $) for even $n$ (odd $n$) and $ \alpha_k$'s are 
arbitrary nonzero constants. The parameters appearing in $ {\bf L} $ 
are $ \delta = 2,\ \gamma = -1 $ and $ \mu = -2. $

  Our final example is the Poincar\'e algebra in $(1+3)$ spacetime. 
The $(1+3)$ Poincar\'e algebra is generated by $P_t$: time translation, 
$ {\bf P} = (P_1, P_2, P_3)$: space translations, $ {\bf J} = (J_1, J_2, J_3)$: 
rotations and $ {\bf K} = (K_1, K_2, K_3)$: boosts. Their commutation 
relations are given by
\beq
\begin{tabular}{lllll}
  $[{\bf J},\; P_t] = 0,$ & $[{\bf J},\; {\bf J}] = {\bf J},$ &
  $[{\bf J},\; {\bf P}] = {\bf P},$ & $[{\bf J},\; {\bf K}] = {\bf K},$ &
  $[P_t,\; {\bf P}] = 0,$  \\
  $[P_t,\; {\bf K}] = {\bf P},$ & $[{\bf P},\; {\bf P}] = 0,$ & 
  $[{\bf K},\; {\bf K}] = -{\bf J},$ & $[{\bf P},\; {\bf K}] = P_t,$& 
\end{tabular}
\label{Pdef}
\eeq 
where $[{\bf P}, {\bf K}] = P_t $ means $ [P_a, P_b] = \delta_{ab} P_t $, 
$[{\bf J}, {\bf J}] = {\bf J}$ means $ [J_a, J_b] = \sum \epsilon_{abc} J_c $ 
with antisymmetric tensor $ \epsilon_{abc} $ and so on. 
The subalgebra $ {\bf L}$ is given by
\bea
 & & J = J_3, \qquad H = K_3, \qquad E  = P_t - P_3, \label{D4P} \\
 & & A = P_1+iP_2, \qquad B = J_1 -iJ_2 + iK_1 + K_2,\nn 
\eea
and other parameters are $ \delta = 1, \gamma = 2i $ and $ \mu = i. $ 
The classical $r$-matrices are obtained from (\ref{cr1})
\beq
 r = (\pm J_3+P_+) \wedge (J_-+iK_-) + 2i K_3 \wedge (P_t-P_3),
 \label{crPoi}
\eeq
where $ J_- \equiv J_1-iJ_2,\ P_+ \equiv P_1+iP_2 $ and $ K_- \equiv K_1 - iK_2.$ 

  Some remarks are in order: (i) One can deform inhomogeneous Lie algebras just 
by twisting their semi-simple or abelian part. 
In the above examples, however, the semi-simple 
part and abelian part are mixed in the each factor of the twisting element $\F.$ 
This indecates that the deformation presented above is nontrivial in the sense 
that the defomation is not caused just by the semi-simple or abelian parts. 
(ii) Our deformation of Poincar\'e algebra is different from the ones in 
\cite{twistedPoin,Mur}. The twisting elements of \cite{twistedPoin} consist of 
only the abelian part, while the twisting element of \cite{Mur} is a product of 
Jordanian twisting factors. (iii) All possible classical $r$-matices for the Poincar\'e 
algebra are classified in \cite{zak}. Our $r$-matrices (\ref{crPoi}) may 
corresponds to the case of $ c= H \wedge X_+,\ a= 0 $ in Table~1 of \cite{zak}.
(iv) We dealt with only Poincar\'e algebras in this article. A classification of 
quantum Poincar\'e grops is given in \cite{PW}.  
(v) The Poincar\'e generators of (\ref{D4P}) are expressed in a 
light-cone type basis. A deformation of the Poincar\'e algebra 
in light-cone type basis (null-plane basis) is discussed 
in \cite{inhom1}. The deformation presented in \cite{inhom1} is simpler 
than the one given by (\ref{D4P}) in the sense that the classical 
$r$-matrix (\ref{crPoi}) has more terms than the one in \cite{inhom1}. 
The peculiarity of the $r$-matirx (\ref{crPoi}) is that it contains 
the term consisiting of rotations : $ J_3 \wedge J_-$. Such term 
is not in the $r$-matirx in \cite{inhom1}.

%
%
%
\section{Conclusion}

  The explicit form of twisting elements that are extensions of PET 
were constructed in this article. Each extensions have two primitive 
elements. It was also shown that our twisting 
elements as well as PET were suitable to deform various inhomogeneous 
Lie algebras. Some of the twisted algebras considered here may admit additional 
Reshetikhin twists \cite{res} that are twisting based on an abelian subalgebra. 
It is shown in \cite{olmo} that additional Reshetikhin twist to PET gives 
an equivalent result as extended twist for the case of $ sl(3).$ It may be 
interesting to consider such mechanism for the twisted inhomogeneous Lie algebras. 
We considered only one parametric deformation. However it is also interesting 
problem to consider possibilities of multiparametric deformation. There 
may exist two different approach to this problem: (i) seeking a possible 
subalgebra (such as $ {\bf L} $ considered in this article) that admits 
multiparametric twisting, (ii) introducing additional twisting elements 
that commute with the ones considered here. The second approach may be 
easier. This will be a future work.

\section*{Acknowledgments}

 The author would like to express his sincere thanks to 
Prof. H.-D. Doebner for his warm hospitality at Technical 
University of Clausthal where this work was done. He also thanks 
Prof. J. Lukierski for helpful comments. This work was supported by 
the Ministry of Education, Science, Sports and Culture, Japan.

\end{document}